\documentclass[10pt]{article}

\usepackage{epsfig,amsmath,amsthm,amssymb}

\usepackage{epsfig,cancel,amsmath,amsthm,amssymb}
\usepackage{color}

\usepackage{lscape}
\usepackage{wrapfig}
\usepackage{color}
\usepackage{graphicx,framed}
\usepackage[colorlinks=true, pdfstartview=FitV, linkcolor=blue, citecolor=blue, urlcolor=blue]{hyperref}
\definecolor{shadecolor}{rgb}{0.95, 0.95, 0.86}

\def\ra{\rightarrow}

\newcommand{\de}{\delta} 
 
\newcommand{\G}{\Gamma} 
\renewcommand{\O}{\Omega}

\newcommand{\e}{\epsilon}
  
\newcommand{\g}{\gamma} 
\newcommand{\la}{\lambda}

\newcommand{\br}{{\mathbb R}}

\newcommand{\CH}{{\mathcal H}}

\newtheorem{theorem}{Theorem}[section]
\newtheorem{example}[theorem]{Example}
\newtheorem{exercise}[theorem]{Exercise}

\newtheorem{lemma}[theorem]{Lemma}
\newtheorem{remark}[theorem]{Remark}

\newtheorem{proposition}[theorem]{Proposition} 
\newtheorem{corollary}[theorem]{Corollary} 
\newtheorem{definition}[theorem]{Definition}
\def\ri{\right}

\def\bth{\begin{theorem}}
\def\et{\end{theorem}}
\def\bc{\begin{corollary}}
\def\ec{\end{corollary}}
\def\bx{\begin{example}\small}
\def\ex{\end{example}}
\def\bxr{\begin{exercise}\small}
\def\exr{\end{exercise}}
\def\bl{\begin{lemma}}
\def\el{\end{lemma}}
\def\bd{\begin{definition}}
\def\ed{\end{definition}}
\def\bp{\begin{proposition}}
\def\ep{\end{proposition}}

\def\br{\begin{remark}}
\def\er{\end{remark}}

\def\be{\begin{equation}}
\def\ee{\end{equation}}
\def\&{\hspace{-15pt}&}
\def\bea{\begin{eqnarray}}
\def\eea{\end{eqnarray}}
\def\beas{\begin{eqnarray*}}
\def\eeas{\end{eqnarray*}}
\def\bi{\begin{itemize}}
\def\ei{\end{itemize}}

\def\C{{\mathbb C}}

\def\R{{\mathbb R}}

\def\wh{\widehat}

\def\a{\alpha}
\def\b{\beta}
\def\e{\epsilon}

\def\l{\lambda}
\def\m{\mu}

\def\1{{\bf 1}}
\def\r{{\rho}}
 
\def\s{ {\sigma}}

\def\th{ {\theta}}

\def\z{\zeta}
\def\x{\xi}

\def\hf{\frac{1}{2}}
\def\le{\left}

\newcommand{\Iscr}{\mathcal I}

\newcommand{\Bscr}{\mathcal B}

\numberwithin{equation}{section}

\begin{document}
\baselineskip 18pt plus 1pt minus 1pt

\begin{center}
\begin{large}

\textbf{\sc Diagonalization of the Finite Hilbert Transform on two adjacent intervals} \\
\end{large}
\bigskip
A. Katsevich$^\star$\footnote{The work  was supported in part by NSF grant DMS-1211164.}
 and 
 A. Tovbis $^\star$\footnote{The work  was supported in part by NSF grant DMS-1211164.}
 \bigskip
 
\begin{small}
$^{\star}$ {\em  University of Central Florida
	Department of Mathematics\\
	4000 Central Florida Blvd.
	P.O. Box 161364
	Orlando, FL 32816-1364
} \\
\end{small}
\end{center}
{\it E-mail:} 
Alexander.Katsevich@ucf.edu, Alexander.Tovbis@ucf.edu \\
\baselineskip 12pt

Abstract. 
We continue the study of stability of solving the interior problem of tomography. The starting point is the Gelfand-Graev formula, which converts the tomographic data into the finite Hilbert transform (FHT) of an unknown function $f$ along a collection of lines. Pick one such line, call it the $x$-axis, and assume that the function to be reconstructed depends on a one-dimensional argument by restricting $f$ to the $x$-axis. Let $\Omega_1$ be the interval where $f$ is supported, and $\Omega_2$ be the interval where the Hilbert transform of $f$ can be computed using the Gelfand-Graev formula.  The equation to be solved is $\left.\CH_1 f=g\right|_{\Omega_2}$, where $\CH_1$ is the FHT that integrates over $\Omega_1$ and gives the result on $\Omega_2$, i.e. $\CH_1: L^2(\Omega_1)\to L^2(\Omega_2)$. In the case of complete data, $\Omega_1\subset\Omega_2$, and the classical FHT inversion formula reconstructs $f$ in a stable fashion. In the case of interior problem (i.e., when the tomographic data are truncated), $\Omega_1$ is no longer a subset of $\Omega_2$, and the inversion problems becomes severely unstable. By using a differential operator $L$ that commutes with $\CH_1$, one can obtain the singular value decomposition of $\CH_1$. Then the rate of decay of singular values of $\CH_1$ is the measure of instability of finding $f$. 

Depending on the available tomographic data, different relative positions of the intervals $\Omega_{1,2}$ are possible. The cases when $\Omega_1$ and $\Omega_2$ are at a positive distance from each other or when they overlap have been investigated already. It was shown that in both cases the spectrum of the operator  $\CH_1^*\CH_1$ is discrete, and the asymptotics of its eigenvalues $\sigma_n$ as $n\to\infty$ has been obtained. In this paper we consider the case when the intervals $\Omega_1=(a_1,0)$ and $\Omega_2=(0,a_2)$ are adjacent. Here $a_1 < 0 < a_2$. Using recent developments in the Titchmarsh-Weyl theory, we show that the operator $L$ corresponding to two touching intervals has only continuous spectrum and obtain two isometric transformations $U_1$, $U_2$, such that $U_2\CH_1 U_1^*$ is the multiplication operator with the function $\sigma(\la)$, $\la\geq(a_1^2+a_2^2)/8$. Here $\la$ is the spectral parameter. Then we show that $\sigma(\la)\to0$ as $\la\to\infty$ exponentially fast. This implies that the problem of finding $f$ is severely ill-posed. We also obtain the leading asymptotic behavior of the kernels involved in the integral operators $U_1$, $U_2$ as $\la\to\infty$. When the intervals are symmetric, i.e. $-a_1=a_2$, the operators $U_1$, $U_2$ are obtained explicitly in terms of hypergeometric functions. \\

{\bf Key words:} interior problem of tomography, finite Hilbert transform, Titchmarsh-Weyl theory, diagonalization\newline

{\bf AMS Subject Classifications:} 34B24, 34L10, 44A12, 44A15.

\section{Introduction}\label{math-intro}

In this paper we continue the study of the stability of solving the interior problem of tomography initiated in papers \cite{kt12, bkt14, aak14, adk15}. The starting point of the study is the Gelfand-Graev formula \cite{gegr-91}, which converts the tomographic data into the finite Hilbert transform (FHT) of an unknown function $f$ along a collection of lines. In what follows we pick one such line, call it the $x$-axis, and assume that the function to be reconstructed depends on a one-dimensional argument by restricting $f$ to the $x$-axis. 

Using the notations of \cite{adk15}, let $\Omega_1$ be the interval where $f$ is supported, and $\Omega_2$ be the interval where the Hilbert transform of $f$ can be computed using the Gelfand-Graev formula.  The equation to be solved can be written in the form $\left.\CH_1 f=g\right|_{\Omega_2}$, where $\CH_1$ is the FHT that integrates over $\Omega_1$ and gives the result on $\Omega_2$, i.e. $\CH_1: L^2(\Omega_1)\to L^2(\Omega_2)$. In the case of complete data, $\Omega_1\subset\Omega_2$, and the classical FHT inversion formula reconstructs $f$ in a stable fashion. In the case of interior problem (i.e., when the tomographic data are truncated), $\Omega_1$ is no longer a subset of $\Omega_2$, and the inversion problem becomes severely unstable. The approach employed in the papers mentioned above is based on a differential operator $L$ that commutes with $\CH_1$. The operator was obtained in \cite{kat10c, kat_11}. By using the commutation property $L \CH_1=\CH_1 L$ one can obtain the singular value decomposition of $\CH_1$. Then the rate of decay of the singular values of $\CH_1$ is the measure of instability of finding $f$. 

Depending on the type of tomographic data available, different relative positions of the intervals $\Omega_{1,2}$ are possible. The case when $\Omega_1$ and $\Omega_2$ are at a positive distance from each other is investigated in \cite{kt12}. It is shown there that the spectrum of the operator  $\CH_1^*\CH_1$ is discrete, and its eigenvalues $\sigma_n$ go to zero exponentially fast as $n\to\infty$. The case when $\Omega_1$ and $\Omega_2$ overlap is investigated in \cite{aak14, adk15}. It is shown that the spectrum of $\CH_1^*\CH_1$ is still discrete and has two accumulation points: 0 and 1. The eigenvalues of the operator can be enumerated in such a way that $\sigma_n\to 0,\, n\to\infty$, and $\sigma_n\to 1,\ n\to -\infty$, and in each case $\sigma_n$ approach the limit exponentially fast. The only case that remained unanswered was when  $\Omega_1$ and $\Omega_2$ touch each other. It was interesting to understand the nature of the spectrum of $\CH_1$ and  estimate how ill-posed it is to find $f$. Since this is a transitional case, it is clear that something special must be happening here. Thus, our problem can be formulated as follows. Given two adjacent intervals $\Omega_1=(a_1,0)$ and $\Omega_2=(0,a_2)$, study the instabily of reconstruction of an $L^2(a_1,0)$ function $f(x)$ knowing its FHT on $(0,a_2)$.

Using recent developments in the Titchmarsh-Weyl theory obtained in \cite{fulton08, be05}, we show in this paper that the operator $L$ corresponding to two touching intervals has only continuous spectrum and obtain two isometric transformations $U_1$, $U_2$, such that $U_2\CH_1 U_1^*$ is a multiplication operator with $\sigma(\la)$, $\la\geq(a_1^2+a_2^2)/8$. Here $\la$ is the spectral parameter. Then we show that $\sigma(\la)\to0$ as $\la\to\infty$ exponentially fast. This implies that the problem of finding $f$ is severely ill-posed. We also obtain the leading asymptotic behavior of the kernels involved in the integral operators $U_1$, $U_2$ as $\la\to\infty$. When the intervals are symmetric, the operators $U_1$, $U_2$ are obtained explicitly in terms of hypergeometric functions. Obviuously, the operator with the kernel $1/(x-y)$ acting from $L^2(-a,0)\to L^2(0,a)$ is naturally related to the operator with the kernel $1/(x+y)$ acting from $L^2(0,a)\to L^2(0,a)$. Thus our results extend those of \cite{rosenblum58}, where, in particular, the diagonalization of the operator $1/(x+y):\, L^2(0,\infty)\to L^2(0,\infty)$ is obtained.  See also the paper \cite{es08}, whether the diagonalization of the operator $1/(x+y):\, L^2(0,\infty)\to L^2(0,\infty)$ is discussed in the context of inverting the Laplace transform.

The paper is organized as follows. In Section~\ref{math-intro} we present the differential operator $L$, establish the commutation relation, and briefly summarize the Titchmarsh-Weyl theory for differential operators with two singular points obtained in \cite{fulton08, be05}. In Section~\ref{sec-Nonsym} we diagonalize the FHT acting from  $L^2(a_1,0)\to L^2(0,a_2)$. In Section~\ref{sym_case} we diagonalize the FHT in the case of symmetric intervals.  In Section~\ref{sec-spectrum} we prove that $L$ does not have discrete spectrum, and some auxiliary results are proven in Sections~\ref{sec-WKB} and \ref{sec-match}.

\section{Spectrum of the commuting differential operator $L$}\label{math-intro}

\subsection{Commuting differential operator}\label{sec-comm}

Fix two points $a_{1,2}$ such that $a_1 < 0 < a_2$ and consider two intervals
\be\label{two-int}
I_1:=(a_1,0),\quad I_2:=(0,a_2).
\ee
Following \cite{kat10c, kat_11}, define a differential operator
\be\label{L-def}
Lf(x)=(P f')'+Qf,\ P(x)=(x-a_1)x^2(x-a_2),\ Q(x)=2\left(x-\frac{a_1+a_2}4\right)^2.
\ee
Each of the intervals gives rise to a singular Sturm-Liouville problem (SLP). By considering the Frobenius solutions near $x=a_1,0$, and $a_2$ we conclude that

\begin{itemize}
\item Near $x=a_j,j=1,2$, there are two linearly independent solutions $\phi_j(x)$ and $\theta_j(x)=\phi_j(x)\ln(x-a_j)+\psi_j(x)$, where $\phi_j$ and $\psi_j$ are analytic near $x=a_j$, and $\psi_j(a_j)=0$;
\item Near $x=0$ there are two linearly independent solutions
\be\label{y-pm}
y_{\pm}=x^{-\hf\pm i\mu}\psi_{\pm}(x)
\ee
where 
\be\label{mu-def}
\mu =\sqrt{ \frac{\l-\frac{(a_1+a_2)^2}8}{-a_1a_2}-\frac 14},
\ee
$\psi_{\pm}(0)=1$, and $\psi_{\pm}(z)$ are analytic in the disk $|z|<\min\{|a_1|,a_2\}$.
\end{itemize}
Also, we see immediately that $L$ is of the Limit Circle (LC) type at $x=a_j,j=1,2$, and of the Limit Point (LP) type at $x=0$. Consequently, no boundary condition is required at $x=0$. The two SLPs become
\be\label{two_slps}
(P(x)f')'+Q(x)f=\lambda f,\ x\in I_j,\ P(x)f'(x)\to 0 \text{ as }x\to a_j,\ j=1,2.
\ee
Next we define two FHTs
\be\label{two_fhts}
(\CH_j f)(z):=\frac1\pi\int_{I_j}\frac{f(x)}{x-z}dx,\ j=1,2.
\ee

\begin{lemma}\label{comm-1} Pick any $f\in C^2(I_j)$, $j=1,2$, such that $f(x)$ is bounded as $x\to a_j$, and
\be\label{f-conds}
f(x)=o(|x|^{-1}),\ f'(x)=o(|x|^{-2}),\ x\to0.
\ee
Then one has:
\be\label{commute}
(\CH_j L f)(x)=(L\CH_j f)(x)\text{ if } \text{dist}(x,I_j)>0,\ j=1,2.
\ee
\end{lemma}

\begin{remark} Assumptions (\ref{f-conds}) are inspired by the properties \eqref{y-pm}, (\ref{mu-def}).
\end{remark}

The proof of the lemma is based on integration by parts and is completely analogous to that of Proposition 2.1 in \cite{kat10c}. The only difference is that now the boundary terms at $x=0$ vanish because \eqref{L-def} and \eqref{f-conds} imply
\be\label{observ}
P(x)f'(x)\to 0 \text{ and } P'(x)f(x)\to 0\text{ as } x\to0.
\ee

%
%

\subsection{Basic facts about diagonalizing the operator $L$\label{diag-basic}}

Consider the operator $L$ acting on smooth functions defined on $I_1$. Recall that $L$ is of the LC type at $a_1$, and of the LP type - at  0. Consider the Liouville transformation 
\be\label{liouv}
t=\int_{a_1}^x\frac{ds}{\sqrt{-P(s)}},\ x\in I_1.
\ee
The transformation \eqref{liouv} maps the interval $I_1$ into the ray $(0,\infty)$. The inverse of the map defines $x=x(t)$ as a function of $t$. Define 
\begin{equation}\begin{split}\label{beta}
F(t):=\sqrt[4]{-P(x)}y(x),\ 
q(t):=Q(x)+\left(\frac{(P'(x))^2}{16P(x)}-\frac{P''(x)}{4}\right),\ x=x(t), t>0.
\end{split}
\end{equation}
A standard computation shows that if $f(x)$ solves the equation $Lf=\lambda f$ on $I_1$, then $F(t)$ solves 
\be\label{transf-eq}
F''(t)+(\lambda-q(t))F(t)=0,\ t>0.
\ee
Note that $q(t(x))\to (a_1^2+a_2^2)/8$ as $x\to 0^-$ (and $t\to\infty$). It is easy to see that after subtracting the constant $(a_1^2+a_2^2)/8$ from $q$ and shifting the spectral parameter accordingly, our potential $q(t)$ satisfies the conditions (1.2)-(1.4) stated in \cite{fulton08}. In particular, in the terminology of \cite{fulton08}, \eqref{transf-eq} falls under Case I with $q_0=1/4$ (cf.  (1.3), (1.4) in \cite{fulton08}). Thus the spectral theory developed in \cite{fulton08}) can be applied to our equation.

Following \cite{fulton08}, we need to find two solutions $\Phi(t,\la),\Theta(t,\la)$ to \eqref{transf-eq} with the following properties:
\begin{equation}\begin{split}\label{props}
&\Phi(t,\la),\Theta(t,\la)\in\mathbb R,\quad \forall t>0,\la\in\mathbb R,\\
&\Phi'(t,\la)\to 0 \text{ as }t\to 0^+,\quad 
W_t(\Theta(t,\la),\Phi(t,\la))=1,\ t>0,\quad \forall\la\in\mathbb C;\\
&\lim_{t\to0} W_t(\Theta(t,\la'),\Phi(t,\la))=1,\quad \forall\la,\la'\in\mathbb C.
\end{split}
\end{equation}
Let  $\phi(x,\la),\theta(x,\la)$ be the solutions to $(L-\la)f=0$ on $I_1$ that correspond to the solutions $\Phi(t,\la),\Theta(t,\la)$ to \eqref{transf-eq}. As is well known, the Wronskians of the two pairs are related by 
\be\label{wronsk}
W_t(\Theta(t,\la),\Phi(t,\la))=-P(x)W_x(\theta(x,\la),\phi(x,\la)).
\ee
Hence, in terms of the solutions to the original equation, conditions \eqref{props} mean:
\begin{equation}\begin{split}\label{props-orig}
&\phi(x,\la),\theta(x,\la)\in\mathbb R,\quad \forall x\in I_1,\la\in\mathbb R,\\
& P(x)\phi'(x,\la)\to 0 \text{ as }x\to a_1^+,\quad
(-P(x))W_x(\theta(x,\la),\phi(x,\la))=1,\ x\in I_1, \quad\forall\la\in\mathbb C;\\
&\lim_{x\to a_1^+} (-P(x))W_x(\theta(x,\la'),\phi(x,\la))=1,\quad \forall\la,\la'\in\mathbb C.
\end{split}
\end{equation}
Note that the first condition on the second line in \eqref{props-orig} is equivalent to the requirement that $\phi(x,\la)$ be bounded as $x\to a_1$ (cf. Lemma~\ref{comm-1}).
Once two solutions $\phi(x,\la),\theta(x,\la)$ that satisfy \eqref{props-orig} have been found, we determine the Titchmarsh-Weyl $m$-function $m(\la)$ from the requirement 
\be\label{m-fn}
\theta(x,\la)+m(\la)\phi(x,\la)\in L^2(I_1),\ \Im\l > 0.
\ee
Then the $m$-function determines the spectral density by the formula
\be\label{rho-density}
\r(\l_2)-\r(\l_1)=\lim_{u\to 0^+} \frac1\pi \int_{\l_1}^{\l_2} \Im m(s+iu)ds,
\ee
where $\l_1,\l_2$ are points of continuity of $\r$. Define the operator $U: L^2(I_1)\to L^2(\mathbb R,d\r)$ and its adjoint by the formulas:
\be\label{u1-def}
(U f)(\l)=\int_{I_1} \phi(x,\l) f(x)dx,\quad (U^* \tilde f)(x)=\int_{\mathbb R} \phi(x,\l) \tilde f(\l)d\r(\l).
\ee
The Titchmarsh-Weyl theory asserts that (cf. \cite{fulton08, be05})
\begin{itemize}
\item the operator $U$ is an isometry: $\Vert f\Vert_{L^2(I_1)}=\Vert U f\Vert_{ L^2(\mathbb R,d\r)}$; 
\item $U$ is unitary: $U^{-1}=U^*$; and 
\item $U$ diagonalizes $L$: ($U L U^{-1}\tilde f)(\l)=\l\tilde f(\l)$ for a sufficiently ``nice" $f$, i.e. for $f\in D(L)$.
\end{itemize}

The interval $I_2$ can be considered in a completely analogous fashion. The only difference is that the two Wronskians in \eqref{props-orig} are multiplied by $P(x)$ instead of $-P(x)$. Thus, the analogue of \eqref{props-orig} becomes
\begin{equation}\begin{split}\label{props-orig-I2}
&\phi(x,\la),\theta(x,\la)\in\mathbb R,\quad \forall x\in I_2,\la\in\mathbb R,\\
& P(x)\phi'(x,\la)\to 0 \text{ as }x\to a_2^-,\quad
P(x)W_x(\theta(x,\la),\phi(x,\la))=1,\ x\in I_2, \quad\forall\la\in\mathbb C;\\
&\lim_{x\to a_2^-} P(x)W_x(\theta(x,\la'),\phi(x,\la))=1,\quad \forall\la,\la'\in\mathbb C.
\end{split}
\end{equation}

\section{General case}\label{sec-Nonsym}

From \eqref{L-def} we have
\be\label{nonsym-ODE-tr1}
x^2(x-a_1)(x-a_2)y''+[2x(x-a_1)(x-a_2)+x^2(2x-a_1-a_2)]y'+[2(x-\frac{a_1+a_2}4)^2-\l] y=0.
\ee

Our first goal is to obtain approximations as $\la\to\infty$ to two linearly independent solutions to \eqref{nonsym-ODE-tr1} that are valid on all $I_1$ and $I_2$. Consider first the interval $I_1=(a_1,0)$. It was shown in \cite{kt12} that in a neighborhood of $x=a_1$ two solutions to \eqref{nonsym-ODE-tr1} can be written in the form
\begin{equation}\label{sols_near_a1}
\begin{split}
& g_1(x)=J_0(2\sqrt{t})+t^{-1/4}O\left(\e^{1-\frac23\de}\right),\ 
g_2(x)=Y_0(2\sqrt{t})+t^{-1/4}O\left(\e^{1-\frac23\de}\right),\\ 
& 1\leq t\leq O\left(\e^{-\left(1+\frac23\de\right)}\right),\ t:=\frac{\l(x-a_1)}{-P'(a_1)},\ \e :=\la^{-1/2},\ 0<\de \ll 1.
\end{split}
\end{equation}
Thus, $t$ is a rescaled variable defined near $x=a_1$. The leading order terms of the WKB solutions to \eqref{nonsym-ODE-tr1}, valid away from $x=a_1,0$, are given by 
\be\label{WKBlead}
\begin{split}
Y_1(x)&=(-P(x))^{-\frac 14}\left\{\cos\le(\sqrt\la \int_{a_1}^x \frac{dt}{\sqrt{-P(t)}}-\frac{\pi}4\ri) + O\left(\e^{\frac12-\de}\right)\right\},\\
Y_2(x)&=(-P(x))^{-\frac 14}\left\{\sin\le(\sqrt\la \int_{a_1}^x \frac{dt}{\sqrt{-P(t)}}-\frac{\pi}4\ri) + O\left(\e^{\frac12-\de}\right)\right\},\\
x&\in [a_1+O\left(\e^{1+2\de}\right),-e^{-1/\sqrt{\e}}].
\end{split}
\ee
Using the asymptotic formulae 8.451.1, 8.451.2 of \cite{GR} for the Bessel functions $J_0(t),Y_0(t)$ as $t\ra\infty$, it was shown in \cite{kt12} that 
\be\label{gY_match}
\begin{split}
g_1(x)&=\frac1{c(\la)}\left\{\left(1+O\left(\e^{\frac12-\de}\right)\right)Y_1(x)+O\left(\e^{\frac12-\de}\right)Y_2(x)\right\},\\
g_2(x)&=\frac1{c(\la)}\left\{O\left(\e^{\frac12-\de}\right)Y_1(x)+\left(1+O\left(\e^{\frac12-\de}\right)\right)Y_2(x)\right\},
\end{split}
\ee
and
\be\label{first_match}
\begin{split}
g_1(x)&=\frac{\cos\le(\varphi(x;\la)-\frac{\pi}4\ri)+O\left(\e^{\frac12-\de}\right)}{c(\la)(-P(x))^{\frac 14}},\
g_2(x)=\frac{\sin\le(\varphi(x;\la)-\frac{\pi}4\ri)+O\left(\e^{\frac12-\de}\right)}{c(\la)(-P(x))^{\frac 14}},\\
x&\in [a_1+O\left(\e^{1+2\de}\right),-e^{-1/\sqrt{\e}}],
\end{split}
\ee
where 
\be\label{c-val}
\varphi(x;\la):=\sqrt \la \int_{a_1}^x \frac{dt}{\sqrt{-P(t)}},\ c(\la):=\la^{1/4}\sqrt{\frac{\pi}{-P'(a_1)}}.
\ee
In a neighborhood of $x=0$ the leading order equation is
\be\label{leading-order}
a_1a_2x^2y''+2a_1a_2xy'+\left[\frac{(a_1+a_2)^2}{8}-\l\right]y=0.
\ee
The characteristic roots are $-\hf\pm i\mu$, where 
\be\label{lam-mu}
\mu =\sqrt{ \frac{\l-\frac{(a_1+a_2)^2}8}{-a_1a_2}-\frac 14}=\sqrt{\frac{\l-\frac{a_1^2+a_2^2}{8}}{-a_1a_2}}=
\sqrt{\frac{\la}{-a_1 a_2}}+O(\e),\ \la\to\infty.
\ee
Thus, $\mu\geq 0$ provided $\l\geq \frac{a_1^2+a_2^2}{8}$.
The corresponding solutions to  \eqref{nonsym-ODE-tr1} have the form 
\be\label{inn-sol}
y_{\pm}(x)=(-x)^{-\hf\pm i\mu}\psi_{1,2}(x;\la),
\ee
where $\psi_{1,2}(0;\la)=1$, and $\psi_{1,2}(x;\la)$ are analytic in the disk $|x|<\max\{|a_1|,a_2\}$.

%
To match the WKB solutions \eqref{WKBlead} with those given in \eqref{inn-sol} we use formula 2.266 of \cite{GR} to obtain
\be\label{eq3}
\begin{split}
\frac1{(-P(x))^{\frac 14}}&=\frac{1+O(x)}{(-a_1 a_2)^{1/4}\sqrt{-x}},\\
\int_{a_1}^x \frac{dt}{\sqrt{-P(t)}}
&=\frac1{\sqrt{-a_1 a_2}}\left.\ln\frac{2\sqrt{-a_1 a_2}\sqrt{(t-a_1)(a_2-t)}+(a_1+a_2)t-2a_1 a_2}{|t|}\right|_{a_1}^{x}\\
&=\frac1{\sqrt{-a_1 a_2}}\left.\ln\frac{[\sqrt{-a_1}\sqrt{a_2-t}+\sqrt{a_2}\sqrt{t-a_1}]^2}{|t|}\right|_{a_1}^{x}\\
&=-\frac{\ln(-x)+\kappa+O(x)}{\sqrt{-a_1 a_2}},\ x\to 0^-,\ 
\kappa:=\ln\frac{a_2-a_1}{-4a_1a_2}.
\end{split}
\ee
Therefore,
\be\label{ph_as}
\varphi(x;\la)=-(\mu+O(\e))(\ln(-x)+\kappa+O(x)).
\ee

For convenience, instead of solutions \eqref{WKBlead} we will temporarily consider an equivalent pair $Y_\pm$:
\be\label{alt-sols}
Y_+(x):=(Y_1(x)+i Y_2(x)) e^{i\pi/4},\ Y_-(x):=(Y_1(x)-i Y_2(x))e^{-i\pi/4}.
\ee
Clearly,
\be\label{WKB-alt}
\begin{split}
Y_\pm(x)&=(-P(x))^{-\frac 14}\left\{\exp(\pm i \varphi(x;\la)) + O\left(\e^{\frac12-\de}\right)\right\},\\
x&\in [a_1+O\left(\e^{1+2\de}\right),-e^{-1/\sqrt{\e}}].
\end{split}
\ee
According to the sentence following \eqref{Voltz} (see Section~\ref{sec-match} below), we will assume $x\in [-c_2\e^2,-c_1\e^2]$, where $0<c_1<c_2<1$. From \eqref{WKBlead}, \eqref{c-val}, and \eqref{ph_as} we find
\be\label{Ypm_as}
\begin{split}
(-x)^{1/2}Y_\pm(x;\la)
&=
\frac{1+O(x)}{(-a_1 a_2)^{1/4} }[\exp(\pm i\varphi(x;\la))+O\left(\e^{\frac12-\de}\right)],\\
&=\frac{1+O(x)}{(-a_1 a_2)^{1/4} }[\exp(\mp i\mu(\ln(-x)+\kappa))+O(\e\ln(-x))+O(x/\e)+O\left(\e^{\frac12-\de}\right)]\\
&=\frac{\exp(\mp i\mu(\ln(-x)+\kappa))+O\left(\e^{\frac12-\de}\right)}{(-a_1 a_2)^{1/4}},\
x\in [-c_2\e^2,-c_1\e^2].
\end{split}
\ee
From \eqref{y-err},
\be\label{ysmall}
(-x)^{1/2}y_\pm(x)=(-x)^{\pm i\mu}+O(\e),\ x\in [-c_2\e^2,-c_1\e^2].
\ee
Matching \eqref{Ypm_as} and \eqref{ysmall} shows
\be\label{match}
Y_\pm(x)=\frac1{(-a_1a_2)^{1/4}}\left\{\exp(\mp i\mu\kappa)y_{\mp}(x)\left(1+O\left(\e^{\frac12-\de}\right)\right)
+\exp(\pm i\mu\kappa)y_{\pm}(x)O\left(\e^{\frac12-\de}\right)\right\}.
\ee
Using \eqref{gY_match}, \eqref{alt-sols}, and \eqref{match} yields
\be\label{eq2}
\begin{split}
g_1(x;\la)= &
\left(\frac{|a_1|^3}{a_2}\right)^{1/4}\sqrt{\frac{a_2-a_1}{\pi}}\frac1{2\la^{1/4}}\biggl\{e^{i\mu\kappa+i\frac\pi4} \left(1+O\left(\e^{\frac12-\de}\right)\right)y_+(x;\la)\\
&\hspace{5cm}+e^{-i\mu\kappa-i\frac\pi4} \left(1+O\left(\e^{\frac12-\de}\right)\right) y_-(x;\la)\biggr\},\\
g_2(x;\la)= &
\left(\frac{|a_1|^3}{a_2}\right)^{1/4}\sqrt{\frac{a_2-a_1}{\pi}}\frac1{2i\la^{1/4}}\biggl\{-e^{i\mu\kappa+i\frac\pi4} \left(1+O\left(\e^{\frac12-\de}\right)\right) y_+(x;\la)\\
&\hspace{5cm}+e^{-i\mu\kappa-i\frac\pi4} \left(1+O\left(\e^{\frac12-\de}\right)\right) y_-(x;\la)\biggr\},\\
\end{split}
\ee


Recall that $t$ is defined according to \eqref{sols_near_a1}. Using formula 8.478 of \cite{GR} it is easy to find that the Wronskian $W_x$ of $\pi Y_0(2\sqrt{t})$ and $J_0(2\sqrt{t})$ (as functions of $x$) 
equals $\frac 1{a_1-x}$. Using the limit
\be
\lim_{x\ra a_1^+} (-P(x)) W_x(\pi Y_0(2\sqrt{t}), J_0(2\sqrt{t}))= -a_1^2(a_2-a_1),
\ee
we obtain that properties \eqref{props-orig} are satisfied if we set
\be\label{phi-theta-1}
\phi_1(x,\la):=g_1(x),\ \theta_1(x,\la):=-\pi g_2(x)/[a_1^2(a_2-a_1)].
\ee
Here and in what follows, the subscript `1' in $\phi_1,\theta_1,m_1$, and $\rho_1$ means that these functions correspond to the interval $I_1$. Condition \eqref{m-fn} now implies that the $m$ function needs to be selected so that the leading coefficients in front of the singularity $(-x)^{-\frac12+i\mu}$ as $x\to 0^-$ in $\theta_1(x,\la)$ and $m_1(\la)\phi_1(x,\la)$ are equal each other in magnitude and are of opposite signs.  Using \eqref{eq2} and \eqref{phi-theta-1} we obtain
\be\label{m-fn-1}
m_1(\la)=\frac{\pi i}{a_1^2(a_2-a_1)}\left(1+O\left(\e^{\frac12-\de}\right)\right).
\ee
Equation \eqref{rho-density} now immediately implies
\be\label{rho-1}
\rho_1'(\l)=\frac1{a_1^2(a_2-a_1)}\left(1+O\left(\e^{\frac12-\de}\right)\right),
\ee
which matches the case of $a_2=-a_1=a$ considered in Section~\ref{sym_case} for large $\l$. 

Next we consider the interval $I_2$. The derivation of the spectral density is very similar, so we sketch here only the main formulas. The analogue of \eqref{sols_near_a1} becomes
\begin{equation}\label{sols_near_a2}
\begin{split}
& g_1(x)=J_0(2\sqrt{t})+t^{-1/4}O\left(\e^{1-\frac23\de}\right),\ 
g_2(x)=Y_0(2\sqrt{t})+t^{-1/4}O\left(\e^{1-\frac23\de}\right),\\ 
& 1\leq t\leq O\left(\e^{-\left(1+\frac23\de\right)}\right),\ t:=\frac{\l(a_2-x)}{P'(a_2)},\ 0<\de \ll 1.
\end{split}
\end{equation}
Thus, $t$ is a rescaled variable defined near $x=a_2$. The leading order terms of the WKB solutions, valid away from $x=0,a_2$, are given by 
\be\label{WKBlead-I2}
\begin{split}
Y_1(x)&=(-P(x))^{-\frac 14}\left\{\cos\le(\sqrt\la \int^{a_2}_x \frac{dt}{\sqrt{-P(t)}}-\frac{\pi}4\ri) + O\left(\e^{\frac12-\de}\right)\right\},\\
Y_2(x)&=(-P(x))^{-\frac 14}\left\{\sin\le(\sqrt\la \int^{a_2}_x \frac{dt}{\sqrt{-P(t)}}-\frac{\pi}4\ri) + O\left(\e^{\frac12-\de}\right)\right\},\\
x&\in [e^{-1/\sqrt{\e}},a_2-O\left(\e^{1+2\de}\right)].
\end{split}
\ee
Matching $g_{1,2}$ in \eqref{sols_near_a2} with $Y_{1,2}$ in \eqref{WKBlead-I2} gives (cf. \cite{kt12}) 
\be\label{gY_match_I2}
\begin{split}
g_1(x)&=\frac1{c(\la)}\left\{\left(1+O\left(\e^{\frac12-\de}\right)\right)Y_1(x)+O\left(\e^{\frac12-\de}\right)Y_2(x)\right\},\\
g_2(x)&=\frac1{c(\la)}\left\{O\left(\e^{\frac12-\de}\right)Y_1(x)+\left(1+O\left(\e^{\frac12-\de}\right)\right)Y_2(x)\right\},
\end{split}
\ee
and
\be\label{first_match_I2}
\begin{split}
g_1(x)&=\frac{\cos\le(\varphi(x;\la)-\frac{\pi}4\ri)+O\left(\e^{\frac12-\de}\right)}{c(\la)(-P(x))^{\frac 14}},\
g_2(x)=\frac{\sin\le(\varphi(x;\la)-\frac{\pi}4\ri)+O\left(\e^{\frac12-\de}\right)}{c(\la)(-P(x))^{\frac 14}},\\
x&\in [e^{-1/\sqrt{\e}},a_2-O\left(\e^{1+2\de}\right)],
\end{split}
\ee
where 
\be\label{c-val-I2}
\varphi(x;\la):=\sqrt \la \int^{a_2}_x \frac{dt}{\sqrt{-P(t)}},\ c(\la):=\la^{1/4}\sqrt{\frac{\pi}{P'(a_2)}}.
\ee
Analogously to \eqref{eq3} we have
\be\label{eq3-I2}
\begin{split}
\int^{a_2}_x \frac{dt}{\sqrt{-P(t)}}
&=-\frac1{\sqrt{-a_1 a_2}}\left.\ln\frac{[\sqrt{-a_1}\sqrt{a_2-t}+\sqrt{a_2}\sqrt{t-a_1}]^2}{|t|}\right|^{a_2}_x\\
&=-\frac{\ln x+\kappa+O(x)}{\sqrt{-a_1 a_2}},\ x\to 0^+,
\end{split}
\ee
where $\kappa$ is the same as in \eqref{eq3}. Therefore,
\be\label{ph_as_I2}
\varphi(x;\la)=-(\mu+O(\la^{-1/2}))(\ln x+\kappa+O(x)).
\ee
With the solutions $Y_\pm$ defined according to \eqref{alt-sols} (using $Y_{1,2}$ for the interval $I_2$), we have
\be\label{WKB-alt-I2}
\begin{split}
Y_\pm(x)&=(-P(x))^{-\frac 14}\left\{\exp(\pm i \varphi(x;\la)) + O\left(\e^{\frac12-\de}\right)\right\},\\
x&\in [e^{-1/\sqrt{\e}},a_2-O\left(\e^{1+2\de}\right)].
\end{split}
\ee
Next we assume $x\in [c_1\e^2,c_2\e^2]$, where $0<c_1<c_2<1$. From \eqref{ph_as_I2} and \eqref{WKB-alt-I2} we find similarly to \eqref{Ypm_as}
\be\label{Ypm_as_I2}
\begin{split}
x^{1/2}Y_\pm(x;\la)
=\frac{\exp(\mp i\mu(\ln x+\kappa))+O\left(\e^{\frac12-\de}\right)}{(-a_1 a_2)^{1/4}},\
x\in [c_1\e^2,c_2\e^2].
\end{split}
\ee
The solutions analogous to \eqref{nonsym-ODE-tr1} have the form 
\be\label{inn-sol-I2}
y_{\pm}(x)=x^{-\hf\pm i\mu}\psi_{1,2}(x;\la),
\ee
where $\psi_{1,2}(0;\la)=1$, and $\psi_{1,2}(x;\la)$ are analytic in the disk $|x|<\max\{|a_1|,a_2\}$. Similarly to \eqref{ysmall},
\be\label{ysmall-I2}
x^{1/2}y_\pm(x)=x^{\pm i\mu}+O(\e),\ x\in [c_1\e^2,c_2\e^2].
\ee
Matching \eqref{Ypm_as_I2} and \eqref{ysmall-I2} shows
\be\label{match-I2}
Y_\pm(x)=\frac1{(-a_1a_2)^{1/4}}\left\{\exp(\mp i\mu\kappa)y_{\mp}(x)\left(1+O\left(\e^{\frac12-\de}\right)\right)
+\exp(\pm i\mu\kappa)y_{\pm}(x)O\left(\e^{\frac12-\de}\right)\right\}.
\ee
Combining \eqref{alt-sols}, \eqref{gY_match_I2}, \eqref{c-val-I2}, and \eqref{match-I2} gives
\be\label{eq2-I2}
\begin{split}
g_1(x;\la)= &
\left(\frac{a_2^3}{|a_1|}\right)^{1/4}\sqrt{\frac{a_2-a_1}{\pi}}\frac1{2\la^{1/4}}\biggl\{e^{i\mu\kappa+i\frac\pi4} \left(1+O\left(\e^{\frac12-\de}\right)\right)y_+(x;\la)\\
&\hspace{5cm}+e^{-i\mu\kappa-i\frac\pi4} \left(1+O\left(\e^{\frac12-\de}\right)\right) y_-(x;\la)\biggr\},\\
g_2(x;\la)= &
\left(\frac{a_2^3}{|a_1|}\right)^{1/4}\sqrt{\frac{a_2-a_1}{\pi}}\frac1{2i\la^{1/4}}\biggl\{-e^{i\mu\kappa+i\frac\pi4} \left(1+O\left(\e^{\frac12-\de}\right)\right) y_+(x;\la)\\
&\hspace{5cm}+e^{-i\mu\kappa-i\frac\pi4} \left(1+O\left(\e^{\frac12-\de}\right)\right) y_-(x;\la)\biggr\},\\
\end{split}
\ee
With $t$ defined according to \eqref{sols_near_a2}, we have $W_x(\pi Y_0(2\sqrt t), J_0(2\sqrt t))=1/(a_2-x)$. Thus 
\be
\lim_{x\ra a_2^-} P(x) W_x(\pi Y_0(2\sqrt t), J_0(2\sqrt t)) = -a_2^2(a_2-a_1),
\ee
and properties \eqref{props-orig-I2} are satisfied by setting
\be\label{phi-theta-I2}
\phi_2(x,\la):=g_1(x),\ \theta_2(x,\la):=-\pi g_2(x)/[a_2^2(a_2-a_1)].
\ee
From \eqref{m-fn}, \eqref{eq2-I2}, and \eqref{phi-theta-I2} we obtain
\be\label{m-fn-I2}
m_2(\la)=\frac{\pi i}{a_2^2(a_2-a_1)}\left(1+O\left(\e^{\frac12-\de}\right)\right).
\ee
Equation \eqref{rho-density} now immediately implies
\be\label{rho-I2}
\rho_2'(\l)= \frac1{a_2^2(a_2-a_1)}\left(1+O\left(\e^{\frac12-\de}\right)\right).
\ee

Now we can find the asymptotics of the diagonal representation of $\CH$. Following \eqref{u1-def} introduce the operators
\be\label{u12-defs}
\begin{split}
(U_1 f)(\l)=\int_{I_1} \phi_1(x,\l) f(x)dx,\quad (U_1^* \tilde f)(x)=\int_{\mathbb R} \phi_1(x,\l) \tilde f(\l)d\r_1(\l),\\
(U_2 f)(\l)=\int_{I_2} \phi_2(x,\l) f(x)dx,\quad (U_2^* \tilde f)(x)=\int_{\mathbb R} \phi_2(x,\l) \tilde f(\l)d\r_2(\l).
\end{split}
\ee
The domain and range spaces of these four operators are defined similarly to Section~\ref{diag-basic}.

Recall that $\phi_{1,2}(x,\l)$ are solutions to $(L-\la)f=0$ on $I_{1,2}$ that are bounded at $a_{1,2}$, respectively. If $\la \geq (a_1^2+a_2^2)/8$, $\phi_{1,2}(x,\l)$ satisfy \eqref{inn-sol-I2}. Thus, $\phi_{1,2}(x,\l)$ satisfy the assumptions of Lemma~\ref{comm-1}, and from \eqref{commute}
\be\label{phi12conn}
\l \CH_1\phi_1=\CH_1 L\phi_1=L\CH_1\phi_1.
\ee
Hence $\CH_1\phi_1$ satisfies on $I_2$ and is bounded near $a_2$. From the Frobenius theory it follows that there cannot be two linearly independent solutions to $(L-\la)f=0$ on $I_2$ that are bounded at $a_2$, so we conclude that $\CH_1\phi_1=\nu(\l)\phi_2$ for some function $\nu(\la)$. Obviously,
\be\label{diag1}
\begin{split}
(\CH_1 U_1^* \tilde f)(\la)=\int_{\mathbb R} (\CH_1\phi_1(x,\l)) \tilde f(\l)d\r_1(\l)
&=\int_{\mathbb R} (\nu(\la)\phi_2(x,\l)) \tilde f(\l)d\r_1(\l)\\
&=\int_{\mathbb R} \left(\nu(\la)\frac{\r'_1(\la)}{\r'_2(\la)}\right)\phi_2(x,\l) \tilde f(\l)d\r_2(\l),
\end{split}
\ee
where $\tilde f\in L^2(\mathbb R,d\r_1)$. In Section~\ref{sec-spectrum} below we will show that $L$ does not have discrete spectrum. It is also well-known that $L$ has no continuous spectrum in the region $\la <(a_1^2+a_2^2)/8$. Hence the integrals in \eqref{diag1} are actually over the interval $\la \geq (a_1^2+a_2^2)/8$. The first equality in \eqref{diag1} holds because $\CH_1:L^2(I_1)\to L^2(I_2)$ is continuous, and the kernel $1/(x-y)$ is smooth on $I_1\times I_2$. Hence
\be\label{diag2}
U_2\CH_1 U_1^* =\nu(\la)\frac{\r'_1(\la)}{\r'_2(\la)},
\ee
To find $\nu(\la)$ we use the well-known identity
\be\label{main-id}
\frac1\pi \int_{-\infty}^0 \frac{(-x)^{-\frac12+i\mu}}{x-y} dx=-\frac{y^{-\frac12+i\mu}}{\cosh(\mu\pi)},\ y>0,\ \mu\in\mathbb R.
\ee
When the interval of integration is not all of $(-\infty,0)$ and the integrand is not exactly $(-x)^{-\frac12+i\mu}$, we can interpret \eqref{main-id} as a statement about the leading singularities. More precisely, if $\CH_1$ acts on a function with the leading singularity $(-x)^{-\frac12+i\mu}$, $x\to 0^-$, the result is a function with the leading singularity $(-1/\cosh(\mu\pi))y^{-\frac12+i\mu}$, $y\to 0^+$. Thus, from \eqref{eq2}, \eqref{phi-theta-1} and \eqref{eq2-I2}, \eqref{phi-theta-I2} we obtain 
\be\label{diag2}
\nu(\la)=-\frac1{\cosh(\mu\pi)}\left(\frac{|a_1|}{a_2^3}\right)^{1/4}\left(\frac{|a_1|^3}{a_2}\right)^{1/4}\left(1+O\left(\e^{\frac12-\de}\right)\right)
=\frac{a_1a_2}{\cosh(\mu\pi)}\left(1+O\left(\e^{\frac12-\de}\right)\right).
\ee
Using now \eqref{rho-1} and \eqref{rho-I2} finally gives
\be\label{diag2}
U_2\CH_1 U_1^* =\frac{a_2^3}{a_1}\frac1{\cosh(\mu\pi)}\left(1+O\left(\e^{\frac12-\de}\right)\right).
\ee

Define $J:=[(a_1^2+a_2^2)/8,\infty)$. The results of this section combined with the results in \cite{fulton08, be05}) can be summarized in the following theorem.

\begin{theorem} The operators $U_j: L^2(I_j)\to L^2(J,\r'_j)$ and $U_j^*: L^2(J,\r'_j)\to L^2(I_j)$, $j=1,2$, defined in \eqref{u12-defs} are isometric transformations. Moreover, in the sense of operator equality on $ L^2(J,\r'_1)$ one has
\be\label{diag_t1}
U_2\CH_1 U_1^* =\sigma(\l),
\ee
where 
\be\label{diag_t2}
\sigma(\l)=\frac{a_2^3}{a_1}\frac1{\cosh(\mu\pi)}\left(1+O\left(\e^{\frac12-\de}\right)\right),\l\to\infty.
\ee
\end{theorem}

\section{Symmetric case}\label{sym_case}

In this section we consider the case of symmetric intervals, i.e. $a_2=-a_1=a$. 
The polynomials $P$ and $Q$ are given by $P=x^2(x^2-a^2)$ and $Q(x)=2x^2$, and the differential equation in \eqref{nonsym-ODE-tr1} becomes
\be\label{sym-ODE-tr1_0}
(x^2(x^2-a^2)y')'+(2x^2-\l)y=0.
\ee
Due to symmetry, if $y(x)$ is a solution to \eqref{sym-ODE-tr1_0}, then so is $y(-x)$.

\subsection{Solution of $Ly=\l y$}\label{sec-sol_L}
The change of variables $x=az$ reduces \eqref{sym-ODE-tr1_0} to 
\be\label{sym-ODE-tr1}
z^2(z^2-1)y''+2z(2z^2-1)y'+ (2z^2-\frac{\l}{a^2})y=0.
\ee
According to \cite{Kamke}, 2.410, two linearly independent solutions of \eqref{sym-ODE-tr1} are given by
\be\label{sol-kamke}
y(z)=z^{-\hf\pm i\mu}\eta_{\pm} (z^2),
\ee
where 
\be\label{alpbetgam}
\mu=\sqrt{\frac{\l}{a^2}-\frac 14}, ~\a=\frac 14 \pm\hf i\mu,~\b=\frac 34 \pm\hf i\mu, ~\g=1\pm i\mu,
\ee
and $\eta_\pm(\xi)$ are solutions of the hypergeometric equation 
\be\label{hyperg}
\xi(\xi-1){\eta''}_\pm+[(\a+\b+1)\xi-\g]\eta'_\pm+\a\b\eta_\pm=0
\ee
with the corresponding choice of the sign in $\a,\b,\g$. Sometimes, we will use notation $\eta$ instead of $\eta_+$.

Since we are interested in a solution $\varphi(z)=\varphi(z,\l)$ of \eqref{sym-ODE-tr1} 
that is analytic at $z=1$, we reduce \eqref{hyperg} to another
hypergeometric equation 
\be\label{hyperg1}
\z(\z-1){\eta''}+[(\a+\b+1)\z-(1+\a+\b-\g)]\eta'+\a\b\eta=0
\ee
by the change of variables $\xi=1-\z$. Then
\be\label{phi}
\varphi(z)=z^{-\hf+i\mu}F(\frac 14 +\frac{ i\mu}2,\frac 34 +\frac{ i\mu}2,1,1-z^2).
\ee
Using the transformation formula 15.3.6 from \cite{AS}, the  behavior of $\varphi$ near $z=0$ is given by
\begin{equation}\label{conn-form}
\begin{split}
\varphi(z)=&\frac {\G(-i\mu)}{\G(\frac 14 -\frac{ i\mu}2)\G(\frac 34 -\frac{ i\mu}2)}
z^{-\hf+i\mu}F(\frac 14 +\frac{ i\mu}2,\frac 34 +\frac{ i\mu}2,1+i \m,z^2)+ \\
&\frac {\G(i\mu)}{\G(\frac 14 +\frac{ i\mu}2)\G(\frac 34 +\frac{ i\mu}2)}
z^{-\hf-i\mu}F(\frac 14 -\frac{ i\mu}2,\frac 34 -\frac{ i\mu}2,1-i \m,z^2)=k f(z) + lg(z),
\end{split}
\end{equation}
where 
\be\label{f-def}
f(z)=z^{-\hf+i\mu}F(\frac 14 +\frac{ i\mu}2,\frac 34 +\frac{ i\mu}2,1+i \m,z^2), \
g(z)=z^{-\hf-i\mu}F(\frac 14 -\frac{ i\mu}2,\frac 34 -\frac{ i\mu}2,1-i \m,z^2),
\ee
and $k,l$ are the corresponding prefactors.

It follows from \eqref{sol-kamke}, \eqref{alpbetgam} that $f(z), ~g(z)$ themselves are solutions to 
\eqref{sym-ODE-tr1} with $f(z)=z^{-\hf+i\mu}\eta_+(z^2)$ and $g(z)=z^{-\hf-i\mu}\eta_-(z^2)$. Moreover, in the case 
\be\label{lamb_ineq}
\l\geq \frac{a^2}4
\ee 
we have $l=\bar k$ and $g(z)=\overline{f(z)}$ when $z\in \R$. Thus, for these values of $\l$ and $z$,
\be\label{real_phi}
\varphi(z,\l)= kf(z)+\bar k \overline{f(\bar z)}=2\Re[kf(z)]. 
\ee
It follows from \eqref{real_phi} that $\varphi(z,\l)$ is real for all $z\in\R$ and $\l\geq a^2/4$. Returning to the original variable $x=az$,
we obtain that
\be\label{phi-x}
\phi(x,\l)=\left(\frac xa\right)^{-\hf+i\mu}F\le(\frac 14 +\frac{ i\mu}2,\frac 34 +\frac{ i\mu}2,1,1-\left(\frac xa\right)^2\ri)
\ee
is a real solution of \eqref{nonsym-ODE-tr1} on $(0,a)$ that is analytic at $x=a$. It is clear that $\phi(-x,\l)$ is also a solution, it is real on $(-a,0)$ and analytic at $x=- a$.

\bl\label{lem-k^2}
If $\l\geq \frac{a^2}{4}$ then 
\be\label{mod_k^2}
|k|^2=\frac{\coth (\pi\mu)}{2\pi\mu}.
\ee
\el

\begin{proof}
 Using \eqref{conn-form}, the Schwarz symmetry of $\G(z)$, and formulae 8.332.1, 8.332.4 of \cite{GR},  we obtain
\be\label{mod_k^2-calc}
|k|^2=\dfrac{|\G(-i\mu)|^2(\frac 14-\frac{i\mu}{2})(\frac 14+\frac{i\mu}{2})}
{\G(\frac 54-\frac{i\mu}{2})\G(\frac 54+\frac{i\mu}{2})\G(\frac 34-\frac{i\mu}{2})\G(\frac 34+\frac{i\mu}{2})}=\frac{\coth (\pi\mu)}{2\pi\mu}.
\ee
\end{proof}

\subsection{Spectral measure for $Ly=\l y$ and diagonalization of $\CH_1$}\label{sec-sol_L}


Following the approach in Section~\ref{sec-Nonsym}, in order to calculate the spectral measure $\r(\l)$ we start with constructing a real-valued solution $\th(x,\l)$, which 
%
must be chosen so that the requirements  \eqref{props-orig} hold. Since  $\th(x)=\th(x,\l)$ must be  linearly independent from $\phi(x,\l)$,
we choose $\th(x,\l)$ as the standard second linearly independent solution of the hypergeometric
equation near $x=a$, see \cite{GR}, 9.153.2,  which can be written as
\be\label{th-sym}
\th(x,\l)=\kappa\left[\phi(x,\l)\ln\left(\frac{a^2-x^2}{a^2}\right)+\Psi\left(\frac{a^2-x^2}{a^2},\l\right)\right],
\ee
where $\Psi(\frac{a^2-x^2}{a^2},\l)$ is the analytic (non-logarithmic) part of this second solution at $x=\pm a$ and $\Psi(0,\l)=0$. We will show below that $\kappa$ is real and $\Psi(\frac{a^2-x^2}{a^2},\l)$ is real-valued for all $x\in\R$ and appropriate $\l$.

\bl\label{lem-th} Set $\kappa=-\frac 1{2a^3}$. 
Let the functions $\phi(x,\l)$ and $\th(x,\l)$ be defined by \eqref{phi-x} and \eqref{th-sym}, respectively. Then the pair $\phi(x,\l)$ and $\th(x,\l)$ satisfies all the requirements \eqref{props-orig} on $(0,a)$, and the pair $\phi(-x,\l)$ and $\th(-x,\l)$ satisfies all the requirements \eqref{props-orig} on $(-a,0)$.
\el

\begin{proof}


We start with the interval $(0,a)$. Clearly, 
\be\label{W1}\begin{split}
W_x:=W_x(\th(x,\l),\phi(x,\l))
&= \kappa \le|
\begin{array}{cc}
\phi\ln\frac{a^2-x^2}{a^2}+\Psi & \phi \\
\phi'\ln\frac{a^2-x^2}{a^2}-\phi\frac{2x}{a^2-x^2}+\frac{d}{dx}\Psi & \phi'
\end{array} \ri|\\
&=\kappa\left(\phi'\Psi-\phi \frac{d}{dx}\Psi+\phi^2\frac{2x}{a^2-x^2}\right). 
\end{split}
\ee
Thus, using that $\phi$ and $\Psi$ are smooth near $a$, we obtain
\be\label{W-th-p}
1=\lim_{x\ra a^-} PW_x(\th,\phi)=-2a^3\kappa.
\ee
Here we have used that $\phi(a,\la)=1$, cf. \eqref{phi-x}. This shows that $\kappa$ is real. By Abel's theorem, $P(x) W_x$ is constant, so the second condition in the second line of \eqref{props-orig} is satisfied. 

Since $\phi$ is real-valued and $P(x)W_x(\theta,\phi)\equiv 1$ on $(0,a)$, the Wronskian of $\phi$ and $\Im \theta$ is zero. Since $\Psi(0,\l)=0$ and $\phi(a,\l)=1$, we immediately conclude that $\Im\Psi\equiv0$.

%
%
Repeating  now the calculations for $W_x(\th(x,\l'),\phi(x,\l))$ and arguing similarly to \eqref{W1}--\eqref{W-th-p},  we obtain
\be\label{W-th-p-d}
\lim_{x\ra a^-}P(x)W(\th(x,\l'),\phi(x,\l))=1
\ee
for any $\l,\l'\in\C$. Note that in this case the logarithmic terms will appear in the Wronskian, but they will not affect the limit in \eqref{W-th-p-d}. Thus our choice of $\kappa$ is correct, and all the requirements in \eqref{props-orig} are satisfied.

Next we consider the interval $(-a,0)$. Analytic continuation of the solutions $\th(x,\l),\phi(x,\l)$, $ \l\geq\frac{a^2}{4}$, from the interval $(0,a)$ to the negative half-axis is no longer real-valued. Therefore, on the interval $(-a,0)$ we replace them by the real-valued solutions $\th(-x,\l),\phi(-x,\l)$. It is straightforward to see that the Wronskian of these solutions is $-\frac{1}{P(-x)}$. However the sign in front of $P(x)$ in \eqref{props-orig} is also changed to the opposite. Thus the pair $\th(-x,\l),\phi(-x,\l)$ satisfies \eqref{props-orig-I2}, and the lemma is proven.
\end{proof}

We are interested in $\Im m(\l)$, where $\l\in\R$.
Given the solutions $\phi$ and $\theta$ with the required properties, we can compute the spectral density $\rho'(\l)$. Again, we start with the interval $(0,a)$. We need $\Im m(\l)$, where $\l\in\R$. In the upper halfplane $\Im \l>0$, the function $m(\l)$ is defined by the requirement that $\th(x,\l)+m(\l)\phi(x,\l)\in L^2(0,a)$, and then $m(\l)$ is analytically continued on the ray $\l\geq\frac{a^2}{4}$. 

Since $\th$ is real-valued and $f(z)$ and $\overline{f(\bar z)}$, where $z=\frac xa$ and $f$ is defined by \eqref{f-def}, are linearly independent, there exists some $l\in\C$ such that  $\th=lf+\bar l \bar f$. 
%
Then,  according to \eqref{f-def} and \cite{AS}, 15.3.10, we have
\be\label{th-ass}
\th(x,\l)=lf(x/a)+\bar l\,\overline {f(x/a)}=-2\Re \left[l\frac {\G(1+i\mu)}{\G(\frac 14 +\frac{ i\mu}2)\G(\frac 34 +\frac{ i\mu}2)}\right]
\ln\frac{a^2-x^2}{a^2}+O(1),~~
x\to a^-.
%
\ee
Note that according to  \eqref{conn-form}, $\frac {\G(1+i\mu)}{\G(\frac 14 +\frac{ i\mu}2)\G(\frac 34 +\frac{ i\mu}2)}=i\mu \bar k$.
Comparing the logarithmic terms of \eqref{th-ass} and
 \eqref{th-sym}, and using Lemma \ref{lem-th},
 we obtain
 \be\label{lbark}
            -[i\mu l\bar k + \overline{i\mu l\bar k}]=-\frac{1}{2a^3}~~~{\rm or}~~~\Im (l\bar k)=-\frac{1}{4\mu a^3}.
 \ee
Let $\Im \l>0$. 
According to \eqref{f-def}, $\bar f \in L^2(0,a)$ and $ f \not\in L^2(0,a)$. So,
the requirement that
\be
\th+m\phi=lf+l\bar f + m(kf+\bar k\bar f)\in L^2(0,a)
\ee
implies $l+mk=0$ or $m=-\frac lk= \frac {-l \bar k}{|k|^2}$. Taking into account \eqref{lbark} and Lemma \ref{lem-k^2},
we obtain
\be\label{m-fn-sym}
\Im m(\l)= \frac {- \Im (l\bar k)}{ |k|^2}=\frac 1{ 4a^3\m |k|^2}=\frac{\pi\tanh (\pi\m)}{2a^3}.
\ee

For the interval $(-a,0)$ and  $\Im \l>0$, the function $m(\l)$ is defined by the requirement that $\th(-x,\l)+m(\l)\phi(-x,\l)\in L^2(-a,0)$. Arguing analogously to \eqref{lbark}--\eqref{m-fn-sym}, we obtain that the $m$-function given in \eqref{m-fn-sym} works for the interval $(-a,0)$ as well. Thus,
\be\label{rho'}
\rho'(\l)=\frac{\tanh (\pi\sqrt{\frac \l{a^2}-\frac 14})}{2a^3},
\ee
and the above holds for both intervals $(-a,0)$ and $(0,a)$.

Using \eqref{conn-form}, we have
\begin{equation}\label{conn-form-conseq}
\begin{split}
\phi_1(x)\sim & k(-x/a)^{-\frac12+i\mu}+\bar k(-x/a)^{-\frac12-i\mu},\ x\to 0^-,\\
\phi_2(x)\sim & k(x/a)^{-\frac12+i\mu}+\bar k(x/a)^{-\frac12-i\mu},\ x\to 0^+.
\end{split}
\end{equation}

Observing that $\rho_1'(\l)/\rho_2'(\l)\equiv1$ (cf. \eqref{diag2}) and combining \eqref{conn-form-conseq} with  \eqref{main-id}, we prove the following result.
\begin {theorem}  Let $J:=[a^2/4)$. Define the functions
\be\label{fns-x-2}
\begin{split}
\phi_2(x,\l)&:=\left(\frac xa\right)^{-\hf+i\mu}F\le(\frac 14 +\frac{ i\mu}2,\frac 34 +\frac{ i\mu}2,1,1-\left(\frac xa\right)^2\ri);\\
\th_2(x,\l)&:=-\frac1{2a^3}\left[\phi_2(x,\l)\ln\left(\frac{a^2-x^2}{a^2}\right)+\Psi\left(\frac{a^2-x^2}{a^2},\l\right)\right],0<x<a,\ \l\in J.
\end{split}
\ee
and 
\be\label{fns-x-1}
\phi_1(x,\l):=\phi_2(-x,\l),\ \th_1(x,\l):=\th_2(-x,\l),\ -a<x<0,\ \l\in J.
\ee
Here $F$ is the hypergeometric function (see 15.1.1 in \cite{AS}), and $\Psi$ is the analytic (non-logarithmic) part of the second solution in \cite{GR}, 9.153.2.
The operators $U_j: L^2(I_j)\to L^2(J,\r')$ and $U_j^*: L^2(J,\r')\to L^2(I_j)$, $j=1,2$, defined in \eqref{u12-defs} are isometric transformations. Moreover, in the sense of operator equality on $ L^2(J,\r')$ one has
\be\label{diag_t2}
U_2\CH_1 U_1^* =\frac{a^2}{\cosh(\mu\pi)}.
\ee
\end{theorem}

\bl\label{lem-k-ass} One has
\be\label{k-ass}
k \sim \frac{e^{i(\frac{\pi}{4}-\mu\ln 2)}}{\sqrt{2\pi\mu}}~~~~~{\rm as}~~~~~\mu\ra +\infty.
\ee
\el

\begin{proof}
The result follows from \eqref{conn-form} and formulae 8.335.1, 8328.2 in \cite{GR}.
\end{proof}

Lemma \ref{lem-k-ass} shows that the behavior of $\phi_2(x;\l)$ as $x\to 0^+$ in the symmetric case (cf. \eqref{conn-form-conseq}) and in the general case (given by \eqref{eq2-I2})
match up.

\subsection{Large $\l$ asymptotics of $\phi(z,\l)$}\label{sec-ass-phi}

In this subsection we calculate a uniform approximation of $\phi(z,\l)$ as $\l\to\infty$. First, we assume for simplicity that $a=1$, so $\varphi(x,\l)=\phi(z,\l)$ and $x=z$. Using \eqref{conn-form} and the integral representation given by formula 9.111 of \cite{GR},  we obtain
\be\label{phi-int}
\phi(z,\l)=\frac{2}{\sqrt{z}}\Re\le[z^{i\m}\frac{\G(1+i\m)\G(-i\m)}{|\G(\frac 14 -\frac{ i\mu}2)\G(\frac 34 -\frac{ i\mu}2)|^2}
\int_0^1e^{\frac{i\m}{2}h(t)}r(t)dt\ri],
\ee
where
\be
h(t)=\ln t+\ln(1-t)-\ln(1-z^2t),~~~r(t)=\frac{1}{[t(1-t)^3(1-z^2t)]^\frac 14}.
\ee
According to Lemma \ref{lem-k^2}, the constant prefactor of the intergral in \eqref{phi-int} is $\frac{i\coth \pi\m}{2\pi}$.
We use the stationary phase method to calculate the asymptotic behavior of the integral.
The stationary point $t_*\in(0,1)$ defined by $h'(t_*)=0$ is calculated to be 
\be
t_*=\frac{1-\sqrt{1-z^2}}{z^2}= \frac{1}{1+\sqrt{1-z^2}}
\ee
We also have 
\be
1-t_*= \frac{\sqrt{1-z^2}}{1+\sqrt{1-z^2}},~~~1-z^2t_*=\sqrt{1-z^2},
\ee
so that 
\be
h(t_*)=-2\ln(1+\sqrt{1-z^2}),~~~     r(t_*)=\frac{1+\sqrt{1-z^2}}{\sqrt{1-z^2}}~~{\rm and}~~~h''(t_*)=-2\frac{(1+\sqrt{1-z^2})^2}{\sqrt{1-z^2}}.
\ee
Applying the stationary phase method and then returning to the original scale (i.e., arbitrary $a$), we get
\be\label{ass-phi}
\phi(x,\l)=\frac{\sqrt{2}a}{\sqrt{\pi \m}\sqrt{x}(a^2-x^2)^{\frac 14}}\cos\le(\m\ln\frac{a+\sqrt{a^2-x^2}}{x}-\frac \pi 4\ri)+O(\mu^{-1}),
\ee
which is valid uniformly on compact subintervals of $(0,a)$. Note that the asymptotics \eqref{ass-phi} in the symmetrical case matches the asymptotics \eqref{first_match_I2} for $\phi_2$ in the general case (cf. \eqref{first_match_I2} and \eqref{phi-theta-I2}). Recall that $\la$ and $\mu$ are related by \eqref{lam-mu}.

\section{Absence of discrete spectrum}\label{sec-spectrum}

In this section we prove that the two Sturm-Liouville problems defined in \eqref{two_slps} have no discrete spectrum. We will consider only the case $j=1$, with the other case being analogous. By assumption, if $\la$ is an eigenvalue and $f(x)$ is the corresponding eigenfunction, then $f$ is bounded (and, hence analytic) near $a_1$ and $f\in L^2(I_1)$. From \eqref{lam-mu} and \eqref{inn-sol} it follows that if $\la >(a_1^2+a_2^2)/8$, then neither of the solutions $y_{\pm}(x)$ is in $L^2(I_1)$. Hence $f\in L^2(I_1)$ imples $f\equiv 0$.  If $\la =(a_1^2+a_2^2)/8$, the solutions behave like $(-x)^{1/2}$ and $(-x)^{1/2}\ln(-x)$, so no linear combination of two such functions can be in $L^2(I_1)$. 

Suppose next that $\la <(a_1^2+a_2^2)/8$. In this case the solutions of $(L-\la)f=0$ behave like $(-x)^{-\frac12\pm q}$ as $x\to 0^-$ for some $q>0$. Clearly, only one of the solutions is in $L^2$. Let $f$ denote the solution which is in $L^2$ and bounded near $a_1$. Thus, $f(x)\sim (-x)^{-\frac12+q}$ as $x\to 0^-$. We can assume $f(a_1)\not=0$, since otherwise $f\equiv0$. Denote $g:=\CH_1 f$. Using \eqref{commute} we have
\be\label{ver1}
\la g=\la \CH_1 f=\CH_1 L f=L\CH_1 f=Lg.
\ee
By the properties of the Hilbert transform, $g$ has the same behavior at zero as $f$: $g(y)\sim y^{-\frac12+q}$ as $y\to 0^+$. Since $Lg=\la g$ on $I_2$, we obtain that $f$ and $g$ are the same solutions up to a multiplicative factor, i.e.
\be\label{one-hilb}
(\CH_1 f)(y)=kf(y),\ y\in I_2,
\ee
where $k$ is a constant. 
Using that $f(a_1)\not=0$ and analytically continuing $f$ from $I_2$ into a neighborhood of $a_1$, we see that $f$ has a logarithmic singularity there. But this contradicts the assumption that $f$ is analytic in a neighborhood of $a_1$. Hence $f\equiv 0$.

\begin{remark} At first glance it follows from equation \eqref{main-id} that $\CH_1$ preserves the ratio of the coefficients in front of the singularities $(-x)^{-\frac12\pm i\mu}$ and, therefore, $\CH_1$ converts a solution of $(L-\la)f=0$ on $I_1$ into (the analytic continuation of) the same solution on $I_2$. This would lead to a contradiction similar to the one obtained above. It is easy to check that $f$ and $\CH_1 f$ are, in fact, two different solutions. Indeed, analytic continuations of $(-x)^{-\frac12\pm i\mu}$ from the negative half-axis to the positive half-axis can be written in the form $c_{\pm}(-x)^{-\frac12\pm i\mu}$, where $c_+\not= c_-$.  Hence the ratios of the coefficients in front of the singularities in $f$ and $\CH_1 f$ at zero are different. 
%
\end{remark}

\section{Validity of the WKB solutions}\label{sec-WKB}

The goal of this section is to construct the WKB solution in a neighborhood of $x=0$.

If equation \eqref{two_slps} is written as a 2 by 2 system, then the 
transformation
\be\label{wkbshear}
Y={\rm diag}(1,\sqrt{\frac{\l}{-P}})\tilde Z
\ee
reduces it to 
\begin{equation}\label{tilZ}
 \e \tilde Z'=\left(
\begin{array}{cc}
0 & \frac{1}{\sqrt{-P}}\\
-\frac{1}{\sqrt{-P}}+\frac{\e^2Q}{\sqrt{-P}}& -\frac{\e P'}{2P}
\end{array}
\right)\tilde Z,
\end{equation} 
where $\e=\frac{1}{\sqrt{\l}}$. Using now 
\begin{equation}
 \tilde Z=\left(
\begin{array}{cc}
1 &i\\
i& 1
\end{array}
\right)Z,
\end{equation} 
we reduce (\ref{tilZ}) to

\begin{equation}\label{matZ}
\e Z'=\left[ \left(
\begin{array}{cc}
 \frac{i}{\sqrt{-P}} &0\\
0& - \frac{i}{\sqrt{-P}}
\end{array}
\right)  -\frac \e 4    \left(
\begin{array}{cc}
\frac{P'}{P}+2i\e \frac Q{\sqrt{-P}} & -i\frac{P'}{P} - 2\e  \frac Q{\sqrt{-P}} \\
i\frac{P'}{P}- 2\e  \frac Q{\sqrt{-P}}  &  \frac{P'}{P}-2i\e \frac Q{\sqrt{-P}} 
\end{array}
\right) \right]Z=AZ.                  
\end{equation} 
Using the Pauli matrices
$$
\s_1= \begin{pmatrix}
0 & 1\\
1&0
\end{pmatrix},~~~
\s_2= \begin{pmatrix}
0 & -i\\
i&0
\end{pmatrix},~~~
\s_3= \begin{pmatrix}
1 & 0\\
0&-1
\end{pmatrix},
$$
we can write $ A=A_0+\e A_1+\e^2A_2$,~~~{\rm where}
\begin{equation}
A_0=\frac{i\s_3}{\sqrt{-P}},~A_1=-\frac{P'}{4P}(I+\s_2),~A_2=\frac{Q}{2\sqrt{-P}}(\s_1-i\s_3).
\end{equation} 
Now the transformation $Z=(I+\e U)X$, where $U=\frac{P'}{8\sqrt{-P}}\s_1$, reduces (\ref{matZ})
to $\e X'=\tilde BX$, where  
\begin{equation}\label{matr-B}
 \tilde B=A_0+\e {\rm diag}A_1+\e^2B(\e)
\end{equation} 
and $B(\e)$ is defined by the equation
\begin{equation}\label{Beps}
(I+\e U) B(\e)=A_2(I+\e U)+A_1U-U{\rm diag}A_1-U'.
\end{equation} 
It is clear that $B(\e)$ is analytic near $\e=0$ provided $\e U$ is small. Direct calculation yields
\begin{equation}\label{B0}
 B(0)=\left(-\frac{(P')^2}{32(-P)^\frac{3}{2}}-\frac Q{2\sqrt{-P}}\right)i\s_3+\left(
\frac{2PP''-(P')^2}{16(-P)^\frac{3}{2}}+\frac Q{2\sqrt{-P}}\right)\s_1.
\end{equation}

Consider equation $\e X'=\tilde BX$ as a perturbation of the diagonal equation
\be
\e W'=(A_0+\e {\rm diag}A_1)W,
\ee
which has a solution 
\be\label{W}
W=P^{-\frac{1}{4}}e^{{\frac{i}{\e}}\int^z\frac{d\z}{\sqrt{-P(\z)}}\s_3}.
\ee
Looking now for a solution of $\e X'=\tilde BX$ in the form $X=TW$, we obtain
\begin{equation}\label{odeT}
 \e T'=\left[A_0+\e {\rm diag}A_1,T\right]+\e^2 BT = \left[A_0,T\right]+\e^2 BT,
\end{equation} 
where 
we have used the fact that ${\rm diag}A_1$ commutes with any matrix $T$ and matrix $W$ is nondegenerate.
Differential equation \eqref{odeT} can be written as Volterra integral equation
\be\label{Voltz}
T(x)=I+\e\int^x e^{\frac {i\s_3}\e \int^x_\z \frac{d\x}{\sqrt{-P(\x)}}}B(\z)T(\z) e^{-\frac {i\s_3}\e \int^x_\z \frac{d\z}{\sqrt{-P(\z)}}}d\z=I+\Iscr T,
\ee
where different contours of integration with the same endpoint $x$  will be selected (see below) for each  entry of the matrix integrand.
We denote this collection of contours by $\tilde \g(x)$. 

\begin{figure}
\begin{center}
\includegraphics[width=0.9\textwidth]{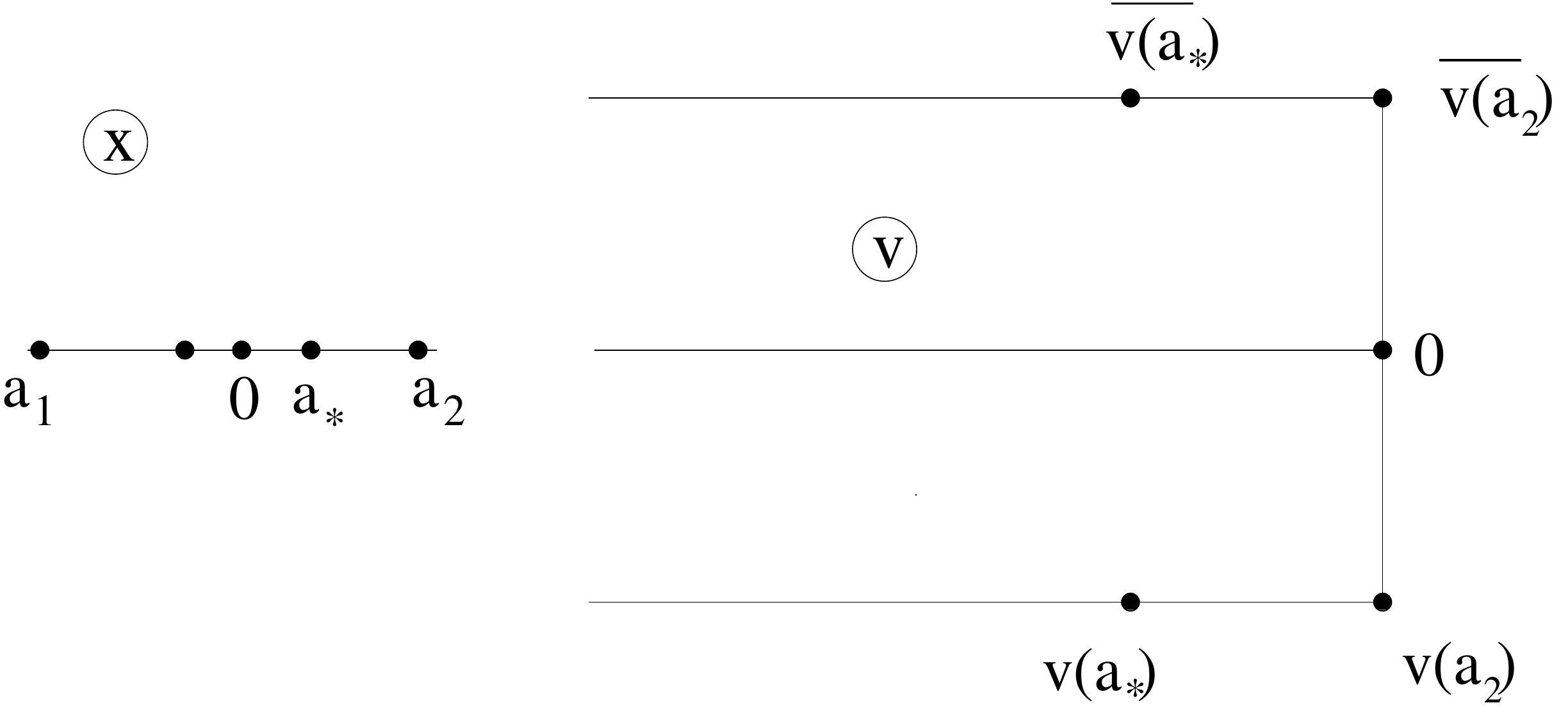}
\caption{The map $v(x)$ maps the complex $x$-plane (left) into the region of the complex $v$-plane, shown on the right. The point
shown on $(a_1,0)$ is $-e^{-\frac{1}{\sqrt{\e}}}$.}
\label{map_v}
\end{center}
\end{figure}

\begin{figure}
\begin{center}
\includegraphics[width=0.9\textwidth]{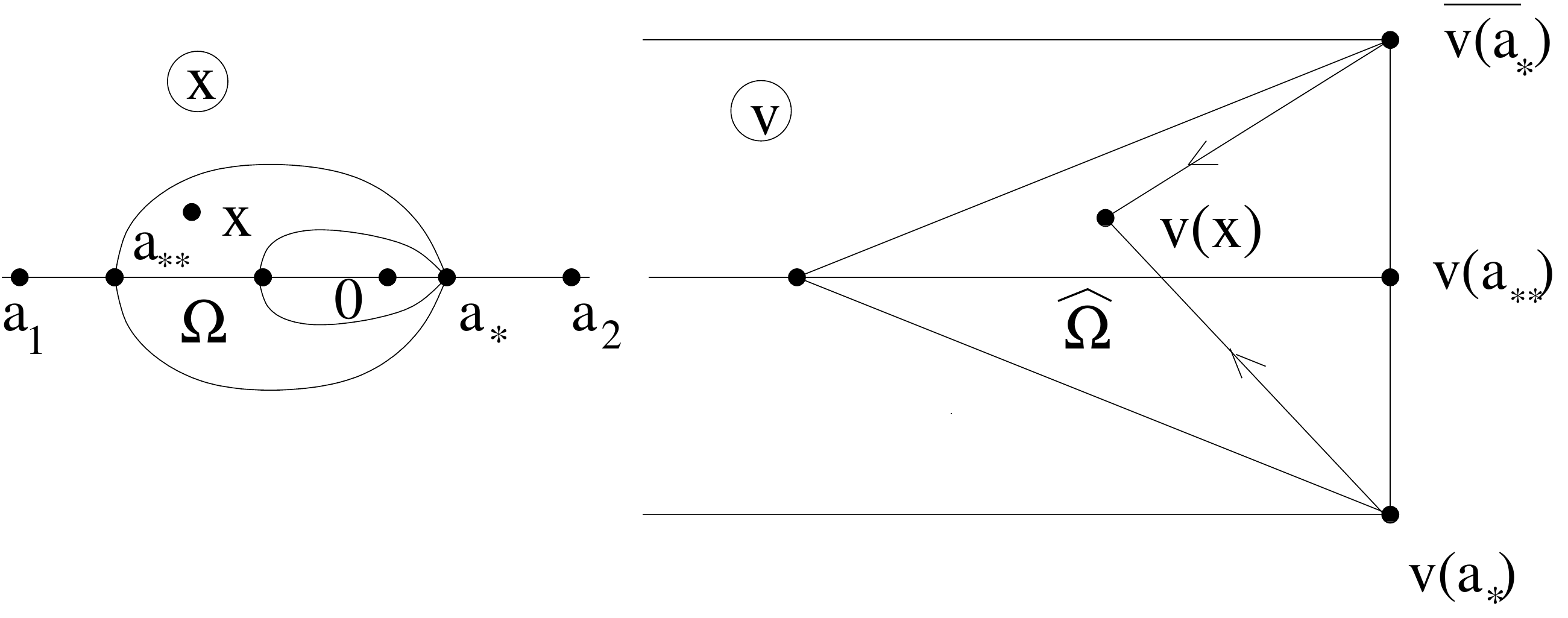}
\caption{The triangular region $\wh\O$.  The preimage $\O$ of $\wh\O$ is shown 
on the left. It has the shape of an oval with a part of its interior (another oval) removed.
Given a point $x\in\O$, the contours $\tilde \g_{1,1}(x)$, $\tilde \g_{2,2}(x)$, are the preimages of the segments  $[v(a_*),v(x)]$, 
$[\overline{v(a_*)},v(x)]$, respectively. The latter are shown on the right. The unmarked points are $-e^{-\frac{1}{\sqrt{\e}}}$ -- on the left, 
and its image $v(-e^{-\frac{1}{\sqrt{\e}}})$ -- on the right.}
\label{Omega}
\end{center}
\end{figure}

We will solve equation \eqref{Voltz} by iterations in a certain region $\O=\O(\e)$ of the complex $x$ plane that comes exponentially close to $x=0$.
In order to  describe the region $\O=\O(\e)$ and contours $\tilde \g(z)$ (and taking into account
\eqref{eq3}), we use the conformal mapping 
\be\label{vmap}
v(x)=\int^x_{a_1}\frac{d\z}{\sqrt{-P(\z)}}
\ee
that maps the upper half plane $\Im x\geq 0$ 
into the semi-strip $-\frac{\pi}{\sqrt{-a_1a_2}}\leq\Im v\leq 0$ and $\Re v\leq 0$ of the complex $v$ plane,
where $v(a_1)=0$, $v(a_2)=-\frac{i\pi}{\sqrt{-a_1a_2}}$ and $v(0)=-\infty$, see Figure \ref{map_v}.
The lower half plane $\Im x\geq 0$  is mapped into the complex conjugated semi-strip. 
Let us pick an arbitrary fixed point $a_*\in(0,a_2)$, for example, $a_*=a_2/2$.
By $\wh\O=\wh\O(\e)$ we define the isosceles triangle with the base $[\overline{v(a_*)},v(a_*)]$ and the (third) vertex at $v(-e^{-\frac{1}{\sqrt{\e}}})$.
According to \eqref{eq3}, 
\be
v(-e^{-\frac{1}{\sqrt{\e}}})=-\e^{-\hf} +O(1) ~~~~~{\rm as}~~~\e\ra 0.
\ee
Then $\O$ is the preimage of $\wh\O$ under the map \eqref{vmap}, which is schematically shown on Figure \ref{Omega}.
It contains the segment $[a_{**},-e^{-\frac{1}{\sqrt{\e}}}]$, $a_{**}\in(a_1,0)$, where $v(a_{**})=\Re v(a_*)$.
Contours $\tilde \g_{1,1}(x)$, $\tilde \g_{2,2}(x)$ are the preimages of the segments $[v(a_*),v(x)]$, 
$[\overline{v(a_*)},v(x)]$. The remaining two
contours connect $\frac{a_1}{2}$ and $x$.

Let $\wh\O_0$, $\O_0$   denote the semi-strip  $|\Im v|\leq \frac{\pi}{\sqrt{-a_1a_2}}$, $\Re v\leq v(a_*)$, and 
its preimage under the map \eqref{vmap}, respectively. Note that $\O_0$ contains both shores of the 
branchcut $[0,a_*]$, and $\O(\e)\subset\O_0$ for all small $\e>0$. Denote by $\Bscr$ the vector space of 
 two by two matrix functions $M(x)$,  which are analytic in $\O_0$ and bounded in $\O(\e)$. The vector space $\Bscr$ becomes a Banach space with the norm given by
$\sup_{x\in\O_0} \Vert M(x)\Vert$, where $\Vert\cdot\Vert$ denotes a matrix norm.

The Volterra equation  \eqref{Voltz} can be written in the operator form as
\be\label{T-iter}
T=({\rm  Id}-\Iscr)^{-1}I=\sum_{j=0}^\infty \Iscr^jI.
\ee
In order to show the convergence of the series in \eqref{T-iter}, we need to estimate the norm of $\Iscr$.
In the variable $v$, the operator $\Iscr$ becomes
\be\label{I-v}
\Iscr M=\e\int_{\wh\g(v)} e^{\frac {i\s_3}\e (v-\x)
}\tilde B(\x) M(\x) e^{-\frac {i\s_3}\e (v-\x)}d\x
\ee
where $\tilde B(\xi)=\left.\sqrt{-P(x)}B(x)\right|_{x=v^{-1}(\xi)}$. According to \eqref{matr-B}-\eqref{B0}, the matrix $\sqrt{-P(x)}B(x)\in\Bscr$.
Let $ \Vert \tilde B(x)\Vert=b$. It follows then from the construction of $\Iscr$ and \eqref{I-v} that
\be\label{Est-n}
\Vert \Iscr M \Vert \leq 2b\e^\hf \Vert M\Vert.
\ee
Thus, choosing $\e<\frac{1}{4b^2}$, we can guarantee the convergence of the series in \eqref{T-iter}, that is, the convergence
of iterations in the solution of the Volterra equation \eqref{Voltz}.

According to the above argument, we have constructed a fundamental solution of the form
\be\label{Y(x)}
Y(x)={\rm diag}\le(1,\sqrt{\frac{\l}{-P}}\ri)\left(
\begin{array}{cc}
1 &i\\
i& 1
\end{array}\right)(I+\frac{\e P'}{8\sqrt{-P}}\s_1)(I+O(\e^{\hf}))(-P)^{-\frac{1}{4}}e^{{\frac{i}{\e}}\int^x_{a_1}\frac{d\z}{\sqrt{-P(\z)}}\s_3}
\ee
on $\O(\e)$. Then, according to \eqref{lam-mu}, \eqref{eq3}, there exist two solution $Y_\pm(x)$ of \eqref{two_slps}, given by
\be\label{Y_pm}
Y_\pm(x)=(-P)^{-\frac{1}{4}}e^{\pm{\frac{i}{\e}}\int^x_{a_1}\frac{d\z}{\sqrt{-P(\z)}}}(I+O(\e^{\hf}))=\frac{e^{\pm i\mu \kappa}}{(-a_1a_2)^\frac 14}(-x)^{-\hf\pm i\mu}(1+O(\e^\hf)+O(x)).
\ee

\section{Validity of the inner solutions}\label{sec-match}

Here we prove the estimate for solutions \eqref{inn-sol}, called inner solutions, on a small interval centered at $x=0$.
This estimate  allows us to match the WKB and inner solutions.

Introducing $y_1=f$, $y_2=Pf'$, we can reduce the original equation \eqref{two_slps} to the matrix equation
\begin{equation}\label{matY}
\tilde Y'=\left(
\begin{array}{cc}
 0 &\frac{1}{P(x)}\\
\l-Q&0
\end{array}
\right) \tilde Y,             
\end{equation} 
where the columns of the matrix $\tilde Y$ are $(y_j,y'_j)$, $j=1,2$, respectively. The shearing transformation
\be\label{innshear}
\tilde Y={\rm diag}(1,\sqrt{\l}x) Y
\ee
reduces \eqref{matY} to
\begin{equation}\label{Y-inn}
\begin{split}
Y'=&\left(
\begin{array}{cc}
0 & \frac{\sqrt \l}{x(x-a_1)(x-a_2)}\\
\frac{\l-Q}{\sqrt{\l}x}& -\frac{1}{x}
\end{array}
\right)Y= (\frac{\tilde B}x+\tilde M)Y\\
=&\left[\frac 1{x} \left(
\begin{array}{cc}
0 &\frac{\sqrt \l}{a_1a_2}\\
\sqrt \l-\frac{(a_1+a_2)^2}{8 \sqrt \l}& - 1
\end{array}
\right) +
\left(
\begin{array}{cc}
0 & \frac{\sqrt {\l}(a_1a+2-x)}{(x-a_1)(x-a_2)}\\
\frac{a_1+a_2-2x}{\sqrt{\l}x}& 0
\end{array}
\right)
\right]Y,
\end{split}
\end{equation} 
where $\tilde B,\tilde M$ are the first and the second terms in the square brackets and $M=M(x)$ is analytic at $x=0$.

It is clear (and can be easily verified) that 
\be\label{Uinn}
\tilde B=U{\rm diag}(-\hf+i\m,-\hf-i\m)U^{-1},~~{\rm where}~~U=
\left(
\begin{array}{cc}
 \frac{1}{a_1a_2} & \frac{1}{a_1a_2}\\
-\frac{1}{2\sqrt \l}+i\frac{\m}{\sqrt\l}&-\frac{1}{2\sqrt \l}-i\frac{\m}{\sqrt\l}
\end{array}
\right),
\ee
and $\m$ is given in \eqref{lam-mu}. The change of variables $Y=UZ$ reduces \eqref{Y-inn} to
$Z'=(\frac{ B}x+ M)Z$, where $B={\rm diag}(-\hf+i\m,-\hf-i\m)$ and $M=U^{-1}\tilde B U$.
Another change of variables $Z=TW$, where $W=x^B$, gives
\be
T'=\frac 1x [i\m \s_3,T]+MT,
\ee
where, according to \eqref{matY}, \eqref{Uinn}, $M=O(\sqrt \l)$. As in Section \ref{sec-WKB},
we replace 
the latter system with the Volterra equation 
\be\label{Voltzi}
T(x)=I+x^{i\m\s_3}\int^x_0 \z^{-i\m\s_3} M(\z)T(\z)\z^{i\m\s_3} d\z x^{i\m\s_3}=I+\Iscr T.
\ee
Since $|x^{\pm i\m}|=1 $ on $\R\setminus\{0\}$, we conclude that on the interval $J=(-\l^{-1},\l^{-1})\subset\R$,
the norm of the operator $\Iscr$ does not exceed $O(\l^{-\hf})$. Thus, we obtain
\be
\tilde Y(x)={\rm diag}(1,\sqrt{\l}x) U(I+O(\l^{-\hf}))x^{-\hf I+i\m\s_3}
\ee
uniformly on $J$. This immediately implies (see \eqref{inn-sol})
\be\label{y-err}
y_{\pm}=(-x)^{-\hf\pm i\mu}(1+O(\l^{-\hf}))
\ee
uniformly on $J$. Since $J$ has a common segment with $\O$ for large $\l$,
we can match the WKB and the inner solutions there. Thus, 
comparing $Y_\pm$ and $y_\pm$ on $\O(\e)$,  we conclude that
\be\label{match_app}
Y_\pm(x)= \frac{e^{\pm i\mu \kappa}}{(-a_1a_2)^\frac 14}y_\pm(x)(1+O(\l^{-\hf})).
\ee


%

\section{Acknowledgement}
AK would like to thank Professor John Schotland, a discussion with whom at the conference ``Mathematical Methods and Algorithms in Tomography" held at the Mathematisches Forschungsinsitut, Oberwolfach, Germany, in August 2014 gave the initial stimulus to this work.



\bibliographystyle{amsalpha}
\bibliography{bibliogr}

\end{document}